# Existence and bifurcation of homoclinic orbits in planar piecewise linear systems

Xiao-Song Yang, Songmei Huan

Department of Mathematics,
Huazhong University of Science and Technology.
Wuhan, 430074, China

**Abstract** The existence and bifurcation of homoclinic orbits in planar piecewise linear homogeneous systems with two regions separated by a discontinuity boundary are investigated in this paper. In addition, existence of periodic orbits and stability of the origin are also discussed without the assumption of continuity for the planar piecewise linear homogeneous systems.

## 1 Introduction

Piecewise-smooth dynamical systems occur naturally in the description of many physical processes, such as impact, friction, switching, sliding and other discrete state transitions and processes with switching components.

There have been a vast number of publications on piecewise-smooth systems in the past several decades. For instance, a lot of books and monographs have been published on this topic. The books [1,2] and [10] dealt with the piecewise-smooth systems from mechanical problems, the detailed analysis of [9] generalized several fundamental theories in smooth system theory including Lyapunov exponent and Conley index to piecewise-smooth cases, and the theories in [12] are related to bifurcations and chaos in the control and electronics systems, and the book [4] gave an comprehensive treatment on the theory and applications of piecewise smooth systems. For journal and conference papers, the reader can refer the recent overview articles [5] for numerous references therein.

Piecewise linear continuous time systems, i.e., piecewise linear ODEs are a typical class of piecewise smooth systems, which are also extensively used to model many physical phenomena. Because piecewise smooth systems are undifferentiable and even discontinuous, the methods from the differentiable dynamical systems theory cannot be applied and specific approaches have to be pursued. Although explicit solutions to a piecewise linear continuous time system can be obtained in each region where one linear system is defined, gluing these solutions globally at the discontinuity boundaries to get a whole picture of the piecewise linear continuous time system is a hard problem, and the theory for piecewise linear continuous time systems is also in its infancy.

The simplest class of piecewise linear continuous time systems of interests should be the planar piecewise linear continuous time systems and will be called planar piecewise linear systems in this paper for brevity. As for the planar piecewise linear systems, several results have appeared in the literature. For example, [6] considered the following so called piecewise linear system (precisely we will call such kind of systems piecewise affine systems)

$$\dot{x} = \begin{cases} A^+ x + b & c^T x > 0 \\ A^- x - b & c^T x \leq 0 \end{cases} \quad (1.1)$$

and provided a complete analysis of uniqueness and non-uniqueness for the initial-value problem, stability and geometrical properties of equilibrium points, existence of sliding motion solutions, periodic orbits and existence of homoclinic orbits. [11] discussed existence of limit cycles and chaotic sets in planar piecewise linear systems under some transition laws.

We will be considering in this paper the simplest case in the sense that the discontinuity boundary is a straight line passing the origin with each component linear system being homogeneous

$$\dot{x} = \begin{cases} A^+ x & c^T x > 0 \\ A^- x & c^T x \leq 0 \end{cases}, \quad x \in R^2 \quad (1.2)$$

where $A^\pm$ be $2 \times 2$ real matrices.

Assume that the overall vector field is continuous across the discontinuity boundary, i.e., $A^+ x = A^- x$, for $x \in \{x : c^T x = 0\}$, then the known facts about system (2) can be summarized as follows.

**Fact 1** [3] Suppose that $(c^T, A^\pm)$ is observable, then the origin is asymptotically stable if and only if

a) neither $A^+$ nor $A^-$ has a real non-negative eigenvalue, and

b) if both $A^+$ and $A^-$ have complex eigenvalues, then $\dfrac{\alpha^+}{\beta^+} + \dfrac{\alpha^-}{\beta^-} < 0$, where $\alpha^\pm \pm \beta^\pm j$, $\beta^\pm > 0$ are eigenvalues of $A^\pm$.

**Fact 2** [3] System (2) has nontrival periodic orbit if and only if both $A^+$ and $A^-$ have complex eigenvalues $\alpha^\pm \pm \beta^\pm j$, $\beta^\pm > 0$ with $\dfrac{\alpha^+}{\beta^+} + \dfrac{\alpha^-}{\beta^-} = 0$, and the period is $\dfrac{\pi}{\beta^+} + \dfrac{\pi}{\beta^-}$.

**Fact 3** [8] If one drops the continuity condition, then Fact 2 holds under the additional condition

$$\langle c, A^+ x \rangle \cdot \langle c, A^- x \rangle > 0 \text{ for } x \in \{x \in R^2 : c^T x = 0\}$$

where $\langle c, A^\pm x \rangle$ is the inner products of $c$ and $A^\pm x$.

Note that the assumption of the observability of $(c^T, A)$ here is equivalent to the assumption that the subspace defined by $c^T x = 0$ is not an invariant submanifold of $\dot{x} = Ax$ in terms of dynamical systems theory.

In view of the above results, it is natural to ask the following questions. Do the above facts hold without assuming that the overall vector field of system (2) is continuous? Is it possible that system (2) has homoclinic orbits for some matrices $A^{\pm}$? With answers to these questions we can provide a scenario of possible dynamics displayed in system (2), and this is the purpose of the present paper. Without loss of generality, we will pay our attentions to the following homogeneous planar piecewise linear system with a discontinuity boundary, written in the form

$$\begin{cases} \dot{x} = A^+ x & x_1 > 0 \\ \dot{x} = A^- x & x_1 \leq 0 \end{cases} \quad x \in R^2 \tag{1.3}$$

here the dots denote derivatives respect to the time $t$, and

$$A^+ = \begin{bmatrix} a_{11}^+ & a_{12}^+ \\ a_{21}^+ & a_{22}^+ \end{bmatrix}, \quad A^- = \begin{bmatrix} a_{11}^- & a_{12}^- \\ a_{21}^- & a_{22}^- \end{bmatrix} \tag{1.4}$$

are non-singular. This is because that if put $c = [c_1 \ c_2]^T$ and define a transformation $T$ as

$$T = \begin{cases} \begin{bmatrix} c_1 & c_2 \\ 0 & sign(c_1) \end{bmatrix}, & c_1 \neq 0 \\ \begin{bmatrix} 0 & c_2 \\ -sign(c_2) & 0 \end{bmatrix}, & c_1 = 0 \end{cases} \tag{1.5}$$

then the change of variables $y = Tx$ transforms system (2) into the form

$$\begin{cases} \dot{y} = B^+ x & y_1 > 0 \\ \dot{y} = B^- x & y_1 \leq 0 \end{cases}$$

where

$$B^{\pm} = TA^{\pm}T^{-1} \text{ and } b_{12} = \begin{cases} \dfrac{1}{c_1 \cdot sign(c_1)} \cdot \det\left(\begin{bmatrix} c^T \\ c^T A \end{bmatrix}\right), & c_1 \neq 0 \\ \dfrac{-1}{c_2 \cdot sign(c_2)} \cdot \det\left(\begin{bmatrix} c^T \\ c^T A \end{bmatrix}\right), & c_1 = 0 \end{cases}.$$

The rest of the paper is organized as follows. In the second section, some preliminaries and the statements of the main results are provided. In the third section, we will prove the main results given in the second section. In the fourth section, the bifurcations of homoclinic orbits are discussed. Finally, some sliding behaviors of system (1.3) are touched on in the last section.

## 2 Preliminaries and Statements of the Main Results

In this section we recall some preliminaries and state the main results. First we recall a notion that is popular in textbooks of control theory.

**Definition 1** Let $C$ be an $m \times n$ matrix, and $A$ be an $n \times n$ matrix, if the rank of the matrix

$$\begin{bmatrix} C \\ CA \\ \vdots \\ CA^{n-1} \end{bmatrix}$$

is $n$, then the pair $(C, A)$ is said to be observable.

**Definition 2** If there exist an equilibrium point $p$ and an orbit $\Gamma$ such that as $t \to \pm\infty$ $\Gamma$ both approaches $p$, then the orbit $\Gamma$ and the equilibrium point $p$ compose a homoclinic orbit.

First we have the following statements concerning the Fact 1.

**Proposition 1** Suppose $(c_0^T, A^\pm)$ is observable, then the origin of system (3) is asymptotically stable if and only if all the following conditions hold.

a) $a_{12}^+ \cdot a_{12}^- > 0$;

b) neither $A^+$ nor $A^-$ has a real non-negative eigenvalue;

c) if both $A^+$ and $A^-$ have complex eigenvalues, then $\dfrac{\alpha^+}{\beta^+} + \dfrac{\alpha^-}{\beta^-} < 0$, where $\alpha^\pm \pm \beta^\pm j$, $\beta^\pm > 0$ are eigenvalues of $A^\pm$.

By virtue of the arguments at the end of the first section, it is easy to get the following result from the above proposition.

**Theorem 1** Suppose $(c^T, A^\pm)$ is observable. If for the overall vector field the conditions b) and c) in the above proposition are satisfied and the following holds

$$\det\left(\begin{bmatrix} c^T \\ c^T A^+ \end{bmatrix}\right) \cdot \det\left(\begin{bmatrix} c^T \\ c^T A^- \end{bmatrix}\right) > 0 \qquad (2.1)$$

then for system (2) the Fact 1 and the Fact 2 above still hold.

**Remarks** i) The condition a) in **Proposition 1** is the special case of (2.1);

ii) It is easy to verify that under the observability condition the continuity of system (2) implies inequality (2.1);

Let $l = \{x \in R^2 : c^T x = 0\}$. In planar case, the inequality (2.1) is in fact equivalent to

$\langle c, A^+ x \rangle \cdot \langle c, A^- x \rangle > 0$ for $x \in l$, which means that the orbits of $\dot{x} = A^+ x$ and $\dot{x} = A^- x$ go from the same one side of $l$ to the other one.

Now we give the main results of this paper.

**Proposition 2** System (3) has homoclinic orbits if and only if the following conditions hold,

a) $(c_0^T, A^\pm)$ is observable, and $A^+, A^-$ both have non-zero real eigenvalues with the same sign;

b) $\lambda^+ \cdot \lambda^- < 0$, $a_{12}^+ \cdot a_{12}^- > 0$, where

$$\lambda^\pm = \begin{cases} \lambda_1^\pm & |\lambda_1^\pm| \geq |\lambda_2^\pm| \\ \lambda_2^\pm & |\lambda_1^\pm| < |\lambda_2^\pm| \end{cases},$$

and $\lambda_{1,2}^\pm$ are the eigenvalues of $A^\pm$.

If a homoclinic orbit exists, then there exists a continuum of homoclinic orbits for system (3).

By this theorem and the arguments at the end of the first section, we can deduce the following result about system (2).

**Theorem 2** The necessary and sufficient condition for system (2) has homoclinic orbits are

a) $A^+$ and $A^-$ both have non-zero real eigenvalues with the same sign, and

$$\det\left(\begin{bmatrix} c^T \\ c^T A^+ \end{bmatrix}\right) \cdot \det\left(\begin{bmatrix} c^T \\ c^T A^- \end{bmatrix}\right) > 0,$$

b) $\lambda^+ \cdot \lambda^- < 0$, where

$$\lambda^\pm = \begin{cases} \lambda_1^\pm & |\lambda_1^\pm| \geq |\lambda_2^\pm| \\ \lambda_2^\pm & |\lambda_1^\pm| < |\lambda_2^\pm| \end{cases},$$

and $\lambda_{1,2}^\pm$ are the eigenvalues of $A^\pm$.

And if exist, there exists a continuum of homoclinic orbits.

## 3 Proofs of the propositions

By definition 1, the following Lemma can easily be shown.

**Lemma 1** If $c_0 = [1\ 0]^T$, then $(c_0^T, A)$ is observable if and only if $a_{12} \neq 0$.

Before giving the proofs of the above propositions, we first show another lemma which plays an important role in proving the main results.

**Lemma 2** If $(c_0^T, A)$ is observable with $\det(A) \neq 0$, $c_0 = [1\ 0]^T$, then the sign of $a_{12}$ determines the directions of trajectories of system $\dot{x} = Ax$ passing through the line $c_0^T x = 0$, i.e., $x_2$-axis. Moreover

1) When $a_{12} > 0$, all trajectories that have intersection points with $x_2$-axis will pass through it clockwise;

2) When $a_{12} < 0$, all trajectories that have intersection points with $x_2$-axis will pass through it counter-clockwise.

The proof is very simple, for reader's convenience we give it as follows.

**Proof** Let $\bar{n}$ be the gradient vector of the line $c_0^T x = 0$, then $\bar{n} = c_0 = [1\ 0]^T$. The tangent vector to the trajectory of system $\dot{x} = Ax$ at any point $(x_1, x_2) \in \{x \in R^2 \mid c_0^T x = 0, x \neq (0,0)\}$ is given by

$$\overline{g} = Ax = \begin{bmatrix} a_{12}x_2 \\ a_{22}x_2 \end{bmatrix}.$$

The inner products of $\overline{n}$ and $\overline{g}$ is

$$\langle \overline{n}, \overline{g} \rangle = c_0^T Ax = a_{12}x_2.$$

It follows that the trajectories intersecting with $c_0^T x = 0$ can pass through it if and only if $\langle \overline{n}, \overline{g} \rangle \neq 0$ for all $x \in \{x \in R^2 \mid c_0^T x = 0\}$, which equivalent to $a_{12} \neq 0$. Furthermore, it is straightforward to get 1) and 2). □

### 3.1 Proof of Proposition 1

By **Lemma 2**, the "only if" part is straightforward, and we show the "if" part as follows.
Let

$$J_0 = \begin{bmatrix} \alpha & -\beta \\ \beta & \alpha \end{bmatrix}, \quad J_1 = \begin{bmatrix} \lambda & 1 \\ 0 & \lambda \end{bmatrix}, \lambda < 0, \quad J_2 = \begin{bmatrix} \lambda_1 & 0 \\ 0 & \lambda_2 \end{bmatrix}, \lambda_1 < \lambda_2 < 0$$

and $J_\pm$ are the Jordan forms of $A^\pm$. Then $J_+, J_- \in \{J_0, J_1, J_2\}$, and we only need to prove that all the trajectories that have intersection points with $x_2$-axis will approach the origin as $t \to +\infty$.

The case with the Jordan forms of $A^\pm$ are both $J_0$ has been discussed in [3]. Here we only give the proof for the case when $J_+ = J_0, J_- = J_1$, and the other cases can be done with the similar way.

Without loss of generality, we put $a_{12}^+ < 0, a_{12}^- < 0$, i.e., the trajectories of $\dot{x} = A^+ x$ and $\dot{x} = A^- x$ pass through $x_2$-axis both counter-clockwise. From [7], for a given initial state $(0, x_{20}), x_{20} \neq 0$ the solution of system $\dot{x} = A^+ x$ is given by

$$\Gamma^+ : \begin{cases} x_1(t) = x_{20} \cdot e^{\alpha t} \cdot \cos(\beta t + \dfrac{\pi}{2}) \\ x_2(t) = x_{20} \cdot e^{\alpha t} \cdot \sin(\beta t + \dfrac{\pi}{2}) \end{cases}$$

where $\alpha \pm \beta j, \beta > 0$ are the eigenvalues of $A^+$. Therefore the solution of the half system $\dot{x} = A^+ x, x_1 > 0$ can be written as

$$\{(x_1 \; x_2) \in \Gamma^+(t) : t \in (0, \pi/\beta), x_{20} < 0\} \tag{3.1}$$

As shown in Fig 1.

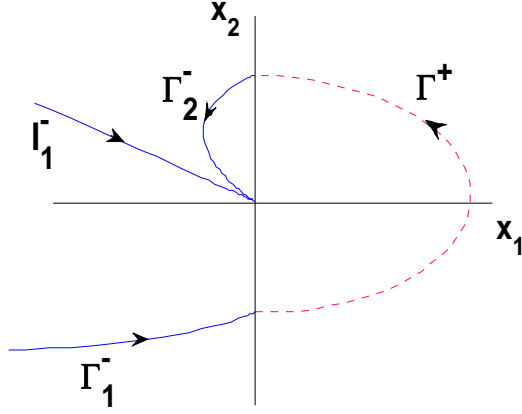

Fig 1   Phase portrait of system (3) with $J_+ = J_0, J_- = J_1$

In order to get the solution of the system $\dot{x} = A^- x$, we let $\xi, \eta$ be the eigenvector and the associated secondary root vector of the eigenvalue $\lambda^-$ of $A^-$. Then $\dim(A^- - \lambda^- I) = 1$, and so there exists nonsingular matrix

$$P = \begin{bmatrix} 1 & 0 \\ p & 1 \end{bmatrix}$$

such that

$$P \cdot (A^- - \lambda^- I) = \begin{bmatrix} m & n \\ 0 & 0 \end{bmatrix},$$

where $m^2 + n^2 \neq 0$. Thus

$$a_{11}^- = \lambda^- + m, \quad a_{12}^- = n, \quad a_{21}^- = -pm, \quad a_{22}^- = \lambda^- - pn \tag{3.2}$$

Submitting $\lambda^- = (a_{11}^- + a_{22}^-)/2$, $(\lambda^-)^2 = \det(A^-)$ into (3.2), we obtain

$$m = pn = \frac{a_{11}^- - a_{22}^-}{2}. \tag{3.3}$$

Moreover

$$(A^- - \lambda^- I)\xi = 0 \Leftrightarrow P(A^- - \lambda^- I)\xi = 0$$
$$(A^- - \lambda^- I)\eta = \xi \Leftrightarrow P(A^- - \lambda^- I)\eta = P\xi,$$

which result in

$$\xi = k \cdot \begin{bmatrix} n \\ -m \end{bmatrix}, \quad \eta = k \cdot \begin{bmatrix} 0 \\ 1 \end{bmatrix}, \quad k \neq 0.$$

Without loss of generality, put $k = 1$. We get by (3.2) and (3.3)

$$M = [\xi \ \eta] = \begin{bmatrix} a_{12}^- & 0 \\ -m & 1 \end{bmatrix} \quad (3.4)$$

Then the transformation $z = M^{-1}x$ changes the system $\dot{x} = A^-x$ into $\dot{z} = J_-z$, and its solution with initial state $(z_{10}, z_{20})$ is given by

$$\begin{cases} z_1(t) = (z_{10} + tz_{20}) \cdot e^{\lambda^- t} \\ z_2(t) = z_{20} \cdot e^{\lambda^- t} \end{cases}$$

Therefore, every solution of system $\dot{x} = A^-x$ with the initial state $(0, x_{20})$, $x_{20} \neq 0$ is

$$\Gamma^- : \begin{cases} x_1(t) = a_{12}^- \cdot x_{20} \cdot t \cdot e^{\lambda^- t} \\ x_2(t) = (x_{20} - m \cdot x_{20} \cdot t) \cdot e^{\lambda^- t} \end{cases} \quad (3.5)$$

And the solution of the half system $\dot{x} = A^-x, x_1 < 0$ can be written as

$$\Gamma^-(t): t \in (-\infty, 0), \quad x_{20} < 0 \quad (3.6)$$

or

$$\Gamma^-(t): t \in (0, +\infty), \quad x_{20} > 0 \quad (3.6)'$$

as shown in Fig 1.

By matching the trajectories given in (3.1), (3.6) and (3.6)', we get the phase portrait of system (3) in Fig 1, where

$$l_1^- : mx_1 + a_{12}^- x_2 = 0 \ (x_1 \leq 0)$$

is the invariant manifold of system (3), and the asymptotical stability of the origin is obvious.

The other case can be discussed in the similar way. For example, with $J_+ = J_2, J_- = J_1$ the phase portrait of system (3) are given in Fig 2

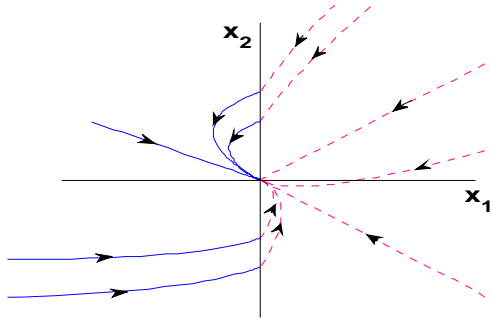

Fig 2   Phase portrait of system (3) with $J_+ = J_2, J_- = J_1$   □

## 3.2 Proof of Proposition 2

Note that by **Definition 2** and the above two lemmas, the " only if" part is obvious and the "if" part can be carried out as follows.

Since the Jordan form of $A^+$ or $A^-$ in question must be one of $J_1$ and $J_2$, where

$$J_1 = \begin{bmatrix} \lambda & 1 \\ 0 & \lambda \end{bmatrix}, \lambda \neq 0, \quad J_2 = \begin{bmatrix} \lambda_1 & 0 \\ 0 & \lambda_2 \end{bmatrix}, \lambda_1 \neq \lambda_2, \lambda_1 \cdot \lambda_2 > 0.$$

we can complete the proof of the "if" part in three cases

**Case 1** Jordan forms of $A^+$, $A^-$ are both $J_1$

For convenience we first consider $\dot{x} = Ax$ with its Jordan form being $J_1$. From the proof of Proposition 1, the solution of the system $\dot{x} = Ax$ with initial state $(x_{10}, x_{20})$ can be written as

$$\begin{cases} x_1(t) = [x_{10} + (mx_{10} + a_{12}x_{20}) \cdot t] \cdot e^{\lambda t} \\ x_2(t) = [x_{20} - \dfrac{m}{a_{12}}(mx_{10} + a_{12}x_{20}) \cdot t] \cdot e^{\lambda t} \end{cases}.$$

where $m = \dfrac{a_{11} - a_{22}}{2}$. According to the signs of $a_{12}, m$ and $\lambda$, the phase portrait of $\dot{x} = Ax$ can easily be sketched in Fig 3. For $M$ given in (3.4), it can be seen that the transformation $z = M^{-1}x$ transforms $z_1$-axis and $z_2$-axis into $l_1 : mx_1 + a_{12}x_2 = 0$ and $l_2 : x_1 = 0$, respectively. The invariant manifold $l_1$ for system $\dot{x} = Ax$ is shown in Fig 3.

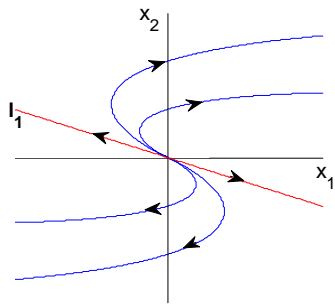 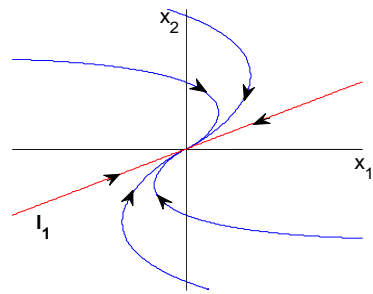

$a_{12} > 0, m > 0, \lambda > 0$ $\qquad\qquad\qquad a_{12} > 0, m < 0, \lambda < 0$

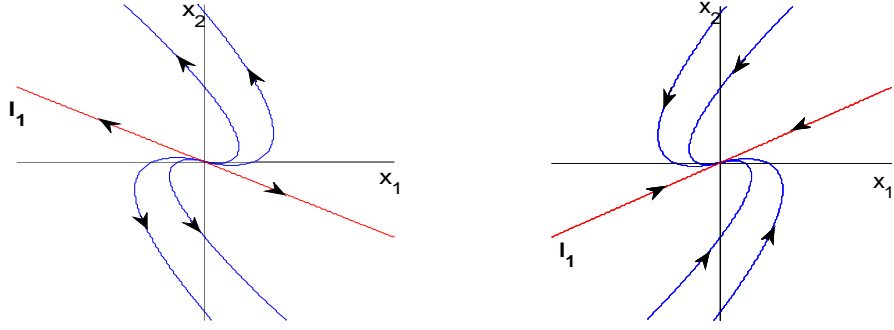

$a_{12} < 0, m < 0, \lambda > 0$  $\qquad\qquad a_{12} < 0, m > 0, \lambda < 0$

Fig 3   Phase portrait of $\dot{x} = Ax$

Now consider (4). Recall that $\lambda^{\pm}$ are the eigenvalues of $A^{\pm}$ and $a_{12}^{\pm}$ are the elements given in (4). For the initial state $p(0, x_{20})$, $x_{20} \neq 0$. By (3.5), the corresponding solutions of system $\dot{x} = A^{-}x$ and system $\dot{x} = A^{+}x$ are given by (see Fig 4)

$$\Gamma_p^{-} : \begin{cases} x_1(t) = a_{12}^{-} \cdot x_{20} \cdot t \cdot e^{\lambda^{-} t} \\ x_2(t) = (x_{20} - m^{-} \cdot x_{20} \cdot t) \cdot e^{\lambda^{-} t} \end{cases} \tag{3.7}$$

$$\Gamma_p^{+} : \begin{cases} x_1(t) = a_{12}^{+} \cdot x_{20} \cdot t \cdot e^{\lambda^{+} t} \\ x_2(t) = (x_{20} - m^{+} \cdot x_{20} \cdot t) \cdot e^{\lambda^{+} t} \end{cases} \tag{3.8}$$

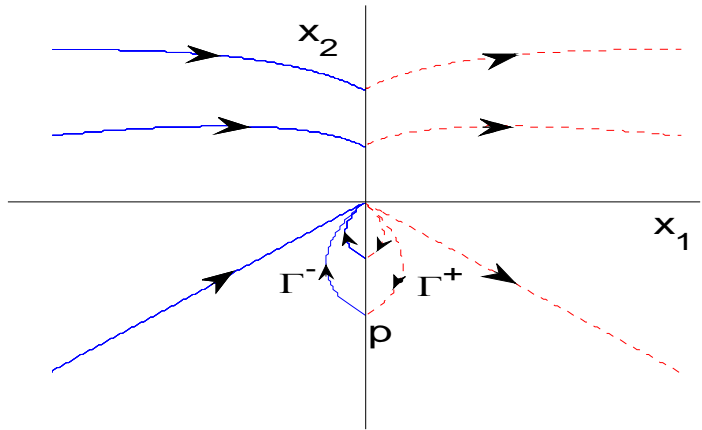

Fig 4   Phase portrait of system (3) with $\lambda^{+} > 0$, $\lambda^{-} < 0$, $a_{12}^{+} > 0$, $a_{12}^{-} > 0$

here $m^{\pm} = \dfrac{a_{11}^{\pm} - a_{22}^{\pm}}{2}$. Furthermore, note that for $\lambda^+ \cdot \lambda^- < 0$ and $a_{12}^+ \cdot a_{12}^- > 0$ there are four cases to be considered.

① $\lambda^+ > 0, \lambda^- < 0$, $a_{12}^+ > 0, a_{12}^- > 0$;   ② $\lambda^+ > 0, \lambda^- < 0$, $a_{12}^+ < 0, a_{12}^- < 0$;

③ $\lambda^+ < 0, \lambda^- > 0$, $a_{12}^+ > 0, a_{12}^- > 0$;   ④ $\lambda^+ < 0, \lambda^- > 0$, $a_{12}^+ < 0, a_{12}^- < 0$.

Without loss of generality, we only prove ①, and the other cases can be done with the same method.

When ① holds, let $x_{20} < 0$, by (3.7) and (3.8), the trajectories of $\Gamma_p^-$ satisfying

$$x_1(t) < 0 \text{ for } t \in (0, +\infty),$$

$$(x_1(t), x_2(t)) \to (0,0) \text{ as } t \to +\infty.$$

and the trajectories of $\Gamma_p^+$ satisfying

$$x_1(t) > 0 \text{ for } t \in (-\infty, 0),$$

$$(x_1(t), x_2(t)) \to (0,0) \text{ as } t \to -\infty.$$

By matching the trajectories $\Gamma_p^+$ and $\Gamma_p^-$, the typical family of trajectories of system (3) are sketched in Fig 4, where

$$l_1^- : m^- x_1 + a_{12}^- x_2 = 0 \ (x_1 < 0), \quad l_1^+ : m^+ x_1 + a_{12}^+ x_2 = 0 \ (x_1 > 0)$$

denote the invariant manifolds of system(3). Thus, the trajectories $\Gamma_p^+(-\infty, 0)$, $\Gamma_p^-(0, +\infty)$ and the point $p$ compose a homoclinic orbit for system (3).

Moreover, it is straightforward to show that if $x(t)$ is a solution of system (1.3), then $k \cdot x(t)$, $k > 0$ is also a solution. Hence, each solution in the set

$$\Sigma = \{(x_1, x_2) \in R^2 \mid m^- x_1 + a_{12}^- x_2 < 0,\ m^+ x_1 + a_{12}^+ x_2 < 0\}$$

is a homoclinic orbit for system (3), and the proof of Case 1 is finished.

**Case 2**  Jordan forms of $A^+, A^-$ are both $J_2$

Putting the Jordan form of $A$ is $J_2$. Let $\lambda_{1,2}$ be its eigenvalues satisfying $|\lambda_1| > |\lambda_2|$, $\lambda_1 \cdot \lambda_2 > 0$, and $\xi_{1,2}$ denote the associated eigenvectors. By simple computing, we get

$$M = [\xi_1 \ \xi_2] = \begin{bmatrix} a_{12} & a_{12} \\ \lambda_1 - a_{11} & \lambda_2 - a_{11} \end{bmatrix}$$

and the transformation $z = M^{-1}x$ changes the system $\dot{x} = Ax$ into $\dot{z} = J_2 z$. Since the the existence of homoclinic orbit is only related to these trajectories that have intersection points with the $x_2$-axis, we give the solution of $\dot{x} = Ax$ with $(0, x_{20})$ be the initial state as follow

$$\begin{cases} x_1(t) = a_{12} \cdot \dfrac{e^{\lambda_2 t} - e^{\lambda_1 t}}{\lambda_2 - \lambda_1} \cdot x_{20} \\ x_2(t) = \dfrac{(\lambda_2 - a_{11})e^{\lambda_2 t} - (\lambda_1 - a_{11})e^{\lambda_1 t}}{\lambda_2 - \lambda_1} \cdot x_{20} \end{cases} \quad (3.9)$$

And according to the signs of $\lambda_1$ and $a_{12}$, the typical trajectories are sketched in Fig 5.

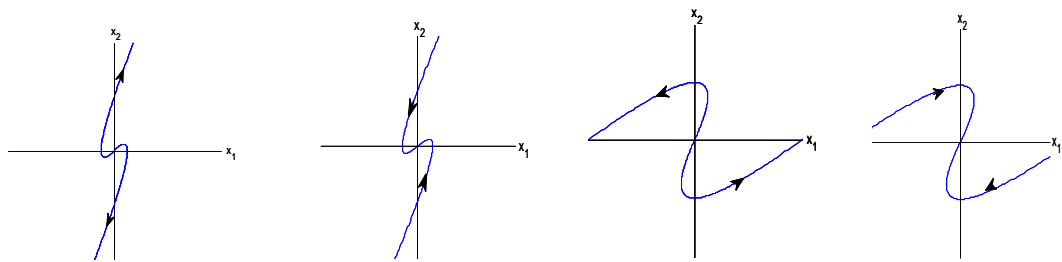

(a) $\lambda_1 > 0, a_{12} > 0$    (b) $\lambda_1 < 0, a_{12} < 0$    (c) $\lambda_1 > 0, a_{12} < 0$    (d) $\lambda_1 < 0, a_{12} > 0$

Fig 5   Phase portrait of $\dot{x} = Ax$

The rest proof of Case 2 is similar to that of Case 1.

**Case 3**   Jordan forms of $A^+, A^-$ are $J_1, J_2$

By the discussions in the foregoing two cases, it is not difficult to obtain the same conclusion in this case. □

## 4 Bifurcation of homoclinic orbits

Now we consider how homoclinic orbits bifurcate when some parameter varies.

Assume that $\lambda_1^+ \cdot \lambda_1^- < 0$, we put

$$A^+ = \begin{bmatrix} \lambda^+ & \mu \\ 0 & \lambda^+ \end{bmatrix}, \quad A^- = \begin{bmatrix} \lambda^- & \mu \\ 0 & \lambda^- \end{bmatrix}$$

where $\mu$ is a parameter. When $\mu$ varies, it is easy to see that the homoclinic orbits of system (1.3) disappear and then appear, which seems to be a new type of bifurcation for homoclinic orbits. In order to make the process clear, we take $\lambda^+ = 1, \lambda^- = -2$ and sketch all the phase portrait of system (1.3) for $\mu = -10, -2, 0, 2, 10$ in Fig 6, from which it can be seen how the homoclinic orbits change as $\mu$ is varing.

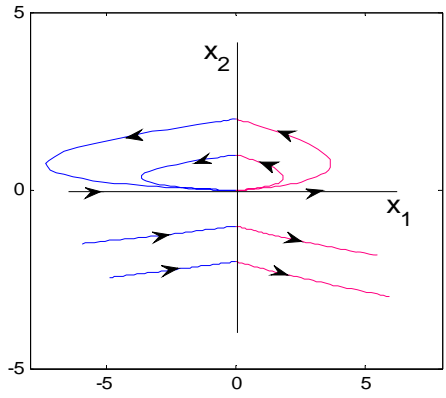
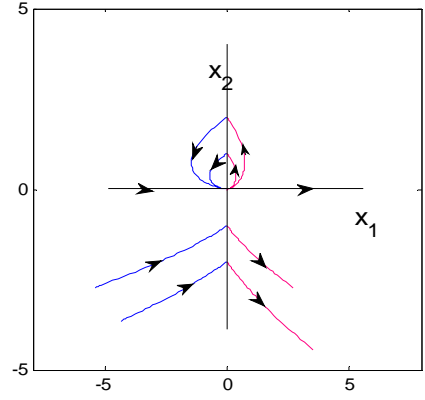

(a) $\mu = -10$        (b) $\mu = -2$

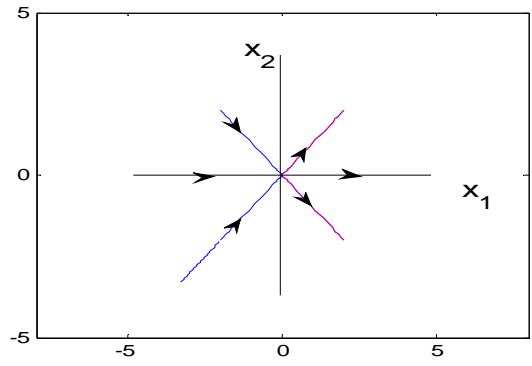

(c) $\mu = 0$

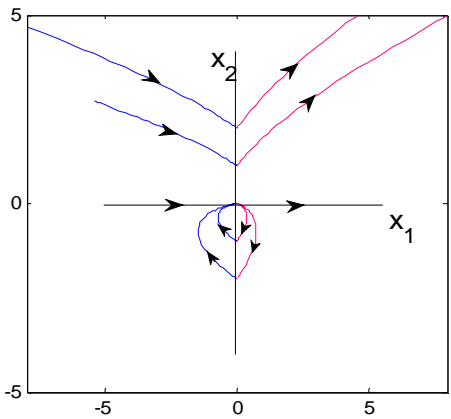
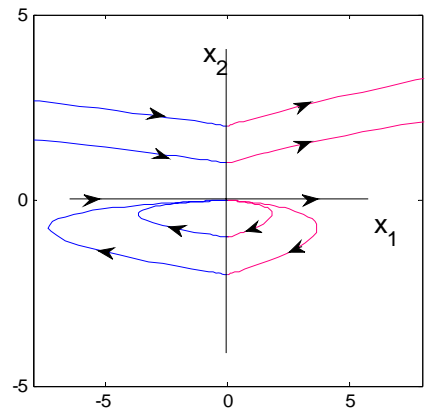

(d) $\mu = 2$        (e) $\mu = 10$

Fig 6    Bifurcation process of system (3)

The above observations can be summarized in the following theorem.

**Theorem 3**    When $\mu < 0$, there always exist a continuum of ***counter-clockwise*** homoclinic orbits above $x_1$-axis (under $x_1$-axis) for system (1.3) with $\lambda^+ > 0(<0)$, $\lambda^- < 0(>0)$, and as $\mu$ is increasing, these homoclinic orbits become more and more "thinner". When $\mu = 0$, all homoclinic orbits disappear, the dynamical behaviors of system (1.3) becomes trivial. When $\mu > 0$, again, there appear a continuum of *clockwise* homoclinic orbits below the $x_1$-axis (above $x_1$-axis) for system (1.3) with $\lambda^+ > 0(<0)$, $\lambda^- < 0(>0)$, and as $\mu$ increasing, these homoclinic orbits become more and more "fatter".

# 5 Conclusions

In this paper we have discussed existence and bifurcation of homoclinic orbits in planar piecewise linear homogeneous systems with two regions separated by a discontinuity boundary. It has been demonstrated in this paper that existence of homoclinic orbits is possible only when $A^+$ and $A^-$ both have nodes. In addition, to expand the results of [3], existence of periodic orbits and asymptotic stability of the origin are also discussed without the assumption of continuity for the planar piecewise linear homogeneous systems. Thus we can have a whole picture about the dynamics of planar piecewise linear homogeneous systems with two pieces separated by a discontinuity boundary. However the situation that the discontinuity boundary contains some sliding regions [4] should be touch upon before completing this paper.

By **Lemma 2**, it is clearly that when $a_{12}^+ \cdot a_{12}^- < 0$, the trajectories of $\dot{x} = A^+ x$ and $\dot{x} = A^- x$ pass through $x_2$-axis in opposite directions, which means the $x_2$-axis is a sliding region for the discontinuity boundary. The sliding region can be "attractive" or "repelling", according to the sign of $\lambda^+ \cdot \lambda^-$, which are summarized in the following paragraph. .

Suppose $(c_0^T, A^\pm)$ is observable, and $A^+, A^-$ both have non-zero real eigenvalues that have the same sign, $a_{12}^\pm$ are the elements of $A^\pm$ defined in (1.4),

$$\lambda^\pm = \begin{cases} \lambda_1^\pm & |\lambda_1^\pm| \geq |\lambda_2^\pm| \\ \lambda_2^\pm & |\lambda_1^\pm| < |\lambda_2^\pm| \end{cases},$$

and $\lambda_{1,2}^\pm$ are the eigenvalues of $A^\pm$. Then we have the following two observations.

**Observation 1** If $a_{12}^+ \cdot a_{12}^- < 0$ and $\lambda^+ \cdot \lambda^- > 0$, the uncertainty of points in $x_2$-axis and the phase portrait of system (1.3) are sketched in Figure 7, where $l_{1,2}^+$ and $l_{1,2}^-$ are the invariant manifolds of system $\dot{x} = A^+ x, x_1 > 0$ and $\dot{x} = A^- x, x_1 \leq 0$, respectively.

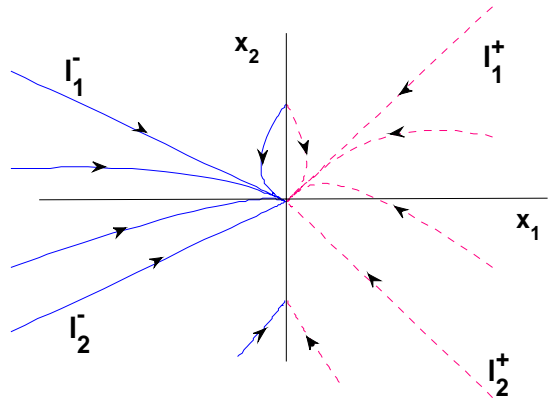

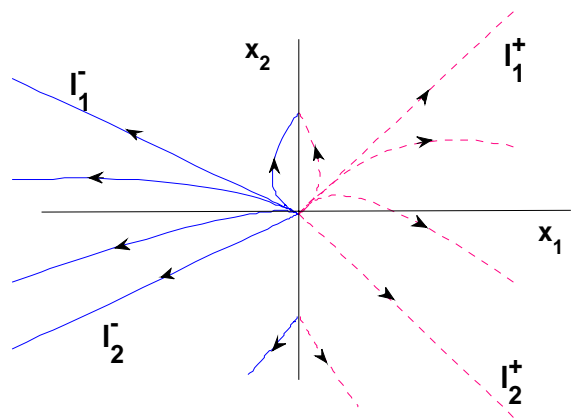

Figure 7

**Observation 2** If $a_{12}^+ \cdot a_{12}^- < 0$ and $\lambda^+ \cdot \lambda^- < 0$, the uncertainty of points in $x_2$-axis and the phase portrait of system (1.3) are sketched in Figure 8, and $l_{1,2}^+$ and $l_{1,2}^-$ are the invariant manifolds of system $\dot{x} = A^+ x, x_1 > 0$ and $\dot{x} = A^- x, x_1 \leq 0$, respectively.

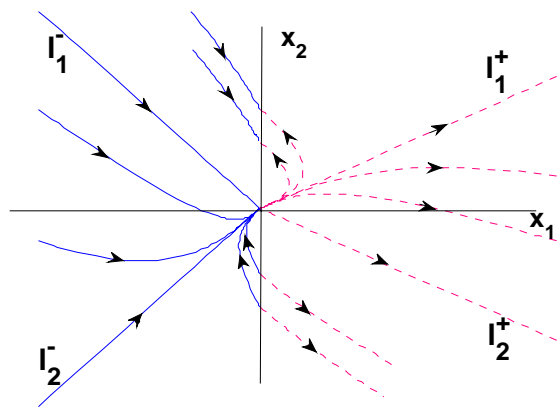

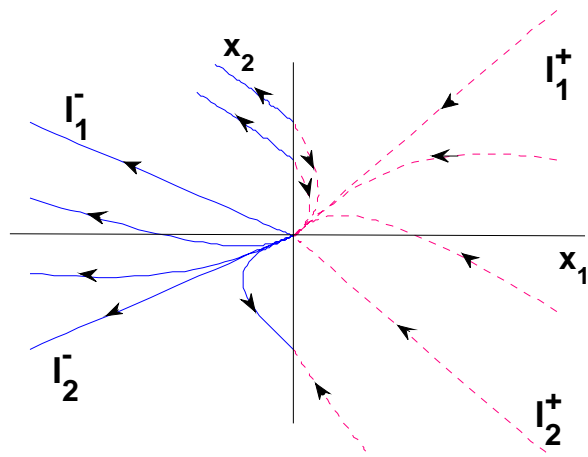

Figure 8

**Acknowledgements** This work is supported in part by National Natural Science Foundation of China （10672062）.